
\hsize=5.7in
\vsize=7.6in
\magnification=\magstep1
\def\proclaim #1. #2\par{\medbreak\vskip-\parskip
    \noindent{\bf#1. \enspace}{\it#2}\par
    \ifdim\lastskip<\medskipamount \removelastskip\penalty55\medskip\fi}

\def\qed {\vrule height4pt width3pt depth2pt}

\def\no{\noindent}

\def\L{{\cal L}}
\def\O{{\cal O}}
\def\N{{\bf N}}
\def\I{{\cal I}}
\def\P{{\bf P}}
\def\Z{{\bf Z}}
\font\srm=cmr8
\centerline{\bf Weierstrass Weight of Gorenstein Singularities with One or Two
Branches}
\bigskip
\centerline{Arnaldo Garcia\footnote{$^1$}{Partially supported by a CNPq grant}
and R.F. Lax\footnote{$^2$}{Partially supported by  NSF grant INT-9101337}}
\medskip
\line{\srm I.M.P.A., Estrada Dona Castorina 110, 22.460 Rio de Janeiro, Brasil
(garcia@impa.br)\hfil}
\line{\srm Department of Mathematics, LSU, Baton Rouge, LA 70803 U.S.A.
(lax@marais.math.lsu.edu)\hfil}
\bigskip
Let $X$ denote an integral, projective Gorenstein curve over an algebraically
closed field $k$.
In the case when $k$ is of characteristic zero, C. Widland and the second
author ([22],
[21], [13]) have defined Weierstrass points of a line bundle on $X$.  In the
first section,
we extend this by defining Weierstrass points of linear systems in arbitrary
characteristic.
This definition may be viewed as a generalization of the definitions of Laksov
[10] and
St\"ohr-Voloch [19] to the Gorenstein case.  Recently  Laksov and Thorup
[11,12] have given
a more general definition of Weierstrass points of ``Wronski systems,''
and our definition may be viewed as a concrete realization in our setting
of their rather abstract definition.

  In the second section, we give an example illustrating our definition.  This
example
  is a plane curve of arithmetic genus 3 in characteristic 2 such that the gap
  sequence at every smooth point (with respect to the dualizing bundle) is
$1,2,5$ and there
  are no smooth Weierstrass points.
    Since every smooth curve of genus 3 in characteristic 2 is classical, this
gives us an
    example of a singular nonclassical curve that is the limit of nonsingular
classical
  curves.  We also compute the Weierstrass weights of points on a rational
curve with a single
  unibranch singularity whose local ring has an especially simple form.

  In the third section, we compute the Weierstrass weight of a unibranch
singularity (on a
  not necessarily rational curve) in terms of its semigroup of values.  In
order to arrive at a
  nice formula, we make the assumption that the characteristic is zero.  We
also compute the
  number of smooth Weierstrass points on a general rational curve with only
unibranch
  singularities.

  In the final section, we compute the Weierstrass weight of a singularity with
precisely two
  branches (again assuming that the characteristic is zero).  This depends
heavily on the
  structure of the semigroup of values of the singularity.  These semigroups
have been studied
  by the first author [4] and by F. Delgado   [1, 2], among others.  In the
course of
  this argument, we construct a basis for the dualizing differentials that is
analogous to a
  (Hermitian) basis of regular differentials adapted to a point in the smooth
case.

\medskip
\no {\bf 1.} Let $k$ denote an algebraically closed field of arbitrary
characteristic.
 Let $X$ be an integral,
projective Gorenstein curve over $k$ of arithmetic genus $g>0$.

 (Gorenstein curves
include any  curve that is  locally a complete intersection; so any curve that
lies on a smooth
surface is Gorenstein.)  Let $K$ denote the field of rational functions on $X$.
Let $\pi:Y\rightarrow X$ denote the normalization of $X$.  Let
$\omega=\omega_X$ denote the
sheaf of dualizing differentials on $X$ and let $\O_{X,P}$ denote the local
ring of the
structure sheaf $\O_X$ at the point $P\in X$.
We recall (cf. [18]) that if $P\in X$, then $\omega_P$ consists of all rational
differentials
$\tau$ on $X$ such that
$$\sum_{Q\rightarrow P}\hbox{Res}_Q(f\tau)=0\quad\hbox{for all }f\in
\O_{X,P},$$
where the sum is over all points on $Y$ lying over $P$.
 Since $X$ is Gorenstein, $\omega$ is an invertible sheaf.

  Let $\L$ be an invertible sheaf on $X$.   Assume that $\L$ has nontrivial
global sections and let
  $V\subseteq  H^0(X,\L)$ be a subspace of dimension $s>0$.
    Choose a
  basis $\phi_1,\phi_2,\dots,\phi_s$ of $V$.  Suppose $\{U_\alpha\}$ is an open
covering
  of $X$ such that $\L(U_\alpha)$ and $\omega(U_\alpha)$ are free
$\O_X(U_\alpha)$-modules
  generated by $\psi_\alpha$ and $\tau_\alpha$, respectively.  Write
  $\phi_j\vert_{U_\alpha} = f_j\psi_\alpha$ for some $f_j\in \O_X(U_\alpha)$
and $j=1,2,\dots,s$.
    Then $f_1,f_2,\dots,f_s$ are linearly independent rational functions over
$k$.

    Suppose $t$ is a separating element of $K$ over $k$ and let $D$ denote the
iterative
    (or Hasse-Schmidt) derivative with respect to $t$.  We recall that the
higher derivatives
    $D^{(i)}$ satisfy the property that
    $$D^{(i)}(\sum_j c_jt^j)=\sum_j {j\choose i} c_j t^{j-i},$$
    where  $c_j\in k$ for all $j$.  F.K. Schmidt [16] showed that there exist
integers
    $$0=\epsilon_0<\epsilon_1<\cdots <\epsilon_{s-1},$$
    minimal with respect to lexicographic ordering of $s$-tuples, such that
    $$\det \pmatrix{D^{(\epsilon_i)} f_j}\ne 0,$$
    where $i=0,1,\dots,s-1$ and $j=1,2,\dots,s$.  These integers are
independent of the choices of
    $\alpha$, the separating element $t$, and the basis of $V$.  The sequence
    $\epsilon_0,\epsilon_1,\dots,\epsilon_{s-1}$ is called the order sequence
of the linear system
     $V$, and one refers to each $\epsilon_i$ as a $V$-order.  If $V$ is a
base-point-free linear
      system, then the order sequence of $V$ is also
     referred to as the order sequence of the morphism
     from $X$ to $\P^{s-1}$ associated to $V$  (cf. [19]), where
$\P^{s-1}=\P_k^{s-1}$ denotes
     the projective space of dimension $s-1$ over $k$.
     If the characteristic of $k$ is 0, then the order sequence is
$0,1,\dots,s-1$.  Put
      $N=\sum_{i=0}^{s-1}\epsilon_i$.

      When the characteristic of $k$ is $p>0$, F. K. Schmidt [16] (also see
[19, Cor. 1.9])
      showed that the order
       sequence of a linear system satisfies the following property.

       \proclaim {(1.1) Proposition}.  Suppose $\epsilon$ is a  $V$-order.  Let
$\mu$ be
       an integer such that
       $${\epsilon\choose \mu}\not\equiv 0\quad(\hbox{mod }p).$$
       Then $\mu$ is also a $V$-order.

       \medskip
       We note that ${\epsilon\choose\mu}\not\equiv 0$ (mod $p$) if and only if
       $\mu\ge 0$ and $\mu$ is $p$-adically smaller than $\epsilon$, which
means that
        each coefficient in the $p$-adic expansion of $\mu$ is less than or
equal to the
	the corresponding coefficient in the $p$-adic expansion of $\epsilon$.
	\medskip
	\no{\bf (1.2) Definition}.  Let $p$ be a prime.
	We say that a finite sequence $\epsilon_0, \epsilon_1,\dots,
	\epsilon_n$ of nonnegative integers satisfies the $p$-adic criterion if
whenever $\mu$ is
	$p$-adically smaller than $\epsilon_i$, for some $i=0,1,\dots,n$, then $\mu$
is a term in the
	sequence (i.e., if the sequence has the property in Proposition (1.1)).
        \proclaim{(1.3) Proposition}.  If $a_0<a_1<\cdots<a_n$ is a sequence
that satisfies
	  the $p$-adic criterion, then these integers are the orders of the morphism
from
	  $\P^1$ to $\P^n$ defined by
	  $$t\mapsto (t^{a_0}:t^{a_1}:\cdots:t^{a_n}).$$

	  \no{\bf Proof}.  [16, Satz 7] (also see M. Homma [8]).\qed

	\medskip

      \def\C{{\cal C}}
      Let
      $\C=\hbox{Ann}(\pi_*{\cal O}_Y/\O_X)=
      \hbox{Ann}(\omega_X/\pi_*\omega_Y)$
       denote the conductor sheaf.
      The dualizing differential $\tau_\alpha$ is of the form
      $$\tau_\alpha= {dt \over h_\alpha},$$ where $t$ is some separating
element for $K$
       over $k$ and where $h_\alpha\in \C(U_\alpha)$.  Put
      $$\rho_\alpha=\det\pmatrix{h_\alpha^{\epsilon_i} D^{(\epsilon_i)} f_j}
      \psi_\alpha^s \tau_\alpha^N,$$
      where $i=0,1,\dots,s-1$ and $j=1,2,\dots,s$.  We note that although
$D^{(\epsilon_i)}f_j$
      may not be in $\O_X(U_\alpha)$, it is in $\pi_*{\cal O}_Y(U_\alpha)$, and
the product
      $h_\alpha^{\epsilon_i}D^{(\epsilon_i)}f_j$ is thus in $\O_X(U_\alpha)$.
Proceeding in
      this way on
      each $U_\alpha$, we obtain functions $\rho_\alpha$ and it is not hard to
show, using
      properties of determinants as in
      [19, Proposition 1.4], that the $\rho_\alpha$ patch to define a section
$\rho\in H^0(X,
      \L^{\otimes s}\otimes \omega^{\otimes N})$, which we refer to as a
wronskian.

      If $P\in X$ and if $\psi$ generates $\L_P$ and $\tau$ generates
$\omega_P$, then we may
      write $\rho=f\, \psi^s\tau^N$ for some nonzero $f\in \O_P=\O_{X,P}$.
        We define ord$_P \rho$ to be
      ord$_Pf=\dim \O_P/(f)$.  This order of vanishing is independent of the
choices of  the
      basis of  $V$ and the generators for $\L_P$ and $\omega_P$.
      \medskip
      \no {\bf (1.4) Definitions.}  Put $W_V(P)=\hbox{ord}_P\rho$ and call this
number the
       $V$-Weierstrass weight of $P$.  The point $P$ is called a
$V$-Weierstrass point if
       $W_V(P)>0$. If $V=H^0(X,\L)$, then we write $W_{\L}(P)$ for  $W_V(P)$
and a
       $V$-Weierstrass point is called an $\L$-Weierstrass point.
       The Weierstrass points of $X$ are the $\omega$-Weierstrass points.  We
will write
       $W_X(P)$, or simply $W(P)$ if it is clear to what curve we are
referring,
        instead of $W_\omega(P)$.
       \medskip
       At any point $P\in X$, one may consider the $(V,P)$-orders, where
       an integer $\mu$ is a $(V,P)$-order if there exists $f\in V$ such that
       ord$_Pf=\mu$.
       At a smooth point $Q$ of $X$, our definition of Weierstrass point
restricts to the
         definitions in  [17], [19],  or [10],   and  one may consider the
        gaps at $Q$, as usual.
	      From [19, Theorem 1.5], we have that if the characteristic of $k$ is
       $p$ and if
        $\epsilon_0(Q),\epsilon_1(Q),\dots,\epsilon_{s-1}(Q)$ are the
       $(V,Q)$-orders, then
       $$W_V(Q)\ge \sum_{i=0}^{s-1}\epsilon_i(Q)-\epsilon_i,$$
       with equality holding if and if
       $$\det \left({\epsilon_j(Q)\choose\epsilon_i}\right)\not\equiv
0\quad(\hbox{mod }p).$$
       (Equality always holds in characteristic 0.)
       \medskip
       \line{{\bf (1.5) Proposition.} {\it The number of $V$-Weierstrass
points, counting multiplicities,
       is}}
       \line{\it $s\deg \L + (2g-2)N$, where $N=\epsilon_0+\epsilon_1+\cdots +
\epsilon_{s-1}$.\hfil}
       \medskip
       \no {\bf Proof.}  This number is the degree of $\L^{\otimes s}\otimes
\omega^{\otimes N}$.\qed
       \medskip
       Let $\O$ denote the local ring at $P$ and let $\tilde\O$ denote the
normalization of $\O$.
       Put $\delta=\delta_P=\dim \tilde\O/\O$.  Suppose that $\tau\in
H^0(X,\omega)$ generates
       $\omega_P$.  Then, locally at $P$, we may write $\tau$  in the form
$\tau=dt/h$, where
       $t$ is a rational function such that ord$_Qt=1$ for all $Q$ on $Y$
       lying over $P$, and where $h$ is some generator (in $\tilde\O$) of the
conductor of
       $\O$ in $\tilde\O$ (cf. [18]).  Since $\O$ is a Gorenstein ring, we have
ord$_Ph=2\delta$.

       \proclaim{(1.6) Proposition}. $W_V(P)=2\delta N +
\hbox{ord}_P\det\pmatrix{D^{(\epsilon_i)}
       f_j}$,  where $i=0,1,\dots,s-1$ and $j=1,2,\dots,s$.

       \no{\bf Proof.} We have
       $$\eqalign{W_V(P)&=\hbox{ord}_P
\det\pmatrix{h^{\epsilon_i}D^{(\epsilon_i)}f_j}\cr
       &=\hbox{ord}_Ph^N +\hbox{ord}_P\det\pmatrix{D^{(\epsilon_i)}f_j}\cr
       &=2\delta N +\hbox{ord}_P\det\pmatrix{D^{(\epsilon_i)}f_j}.\qed\cr}$$

       \proclaim{(1.7) Corollary}. $W_V(P)\ge 2\delta N$. In particular, if $P$
is a singular point and
       if $s>1$, then $P$ is a $V$-Weierstrass point.

       \no {\bf Proof.} If $s>1$, then $N=\epsilon_0+\epsilon_1+\cdots
+\epsilon_{s-1}\ge 0+1=1$.
         Hence $W_V(P)\ge 2\delta>0$, since $P$ is singular.\qed
	 \medskip
	 \goodbreak
	 \no{\bf 2.} As one might expect, phenomena can occur on
	  Gorenstein curves in positive
	 characteristic that do not occur on Gorenstein curves in characteristic zero.
 And  there
	 exist (singular) Gorenstein curves in positive characteristic that exhibit
behavior that is
	 not found on smooth curves.  The following example serves to illustrate both
of these
	 points.
	 \medskip
	 \no {\bf (2.1) Example.} Suppose the characteristic of $k$ is 2.
	  Let $X$ be the rational Gorenstein curve obtained from the
	  projective line over
	 $k$ by replacing the local ring $\tilde \O$ at 0 with the ring
	 $$\O=k+kt^3+kt^4+t^6\tilde \O,$$
	 where $t$ is a uniformizing parameter of $\tilde \O$ that generates the
function
	 field of $X$.  Equivalently, $X$ is the plane quartic
	 $y^3z=x^4$.   Let $P$ denote the singular point of $X$.

	 We will find the Weierstrass points of $X$ (i.e., the $\omega$-Weierstrass
points).
	   It is easy to see that the rational  differentials
	 $$\tau_1=dt/t^6, \tau_2=dt/t^3, \hbox{ and }\tau_3=dt/t^2$$
	 form  a basis of $H^0(X,\omega)$  and that $\tau_1$ generates
	 $\omega_P$.  Referring to the notation in \S 1, we have $f_1=1, f_2=t^3,
f_3=t^4$, and
	 $h=t^6$.
	 Let $(D^{(i)}f_j)$ denote the triple whose components are the $i$-th
iterative derivatives
	 (with respect to $t$) of the functions $f_1,f_2$, and $f_3$.  Then we have
         $$\eqalign{&(D^{(1)}f_j)=(0,t^2,0)\cr &(D^{(2)}f_j)=(0,t,0)\cr
	                   & (D^{(3)}f_j)=(0,1,0)\cr &(D^{(4)}f_j)=(0,0,1)\cr}$$
         Therefore, the order sequence of $\omega$ is 0,1,4.  The wronskian at
$P$ is then
	 $$\det\pmatrix{1&t^3&t^4\cr 0&t^8&0\cr 0&0&t^{24}\cr}=t^{32}.$$
	 It follows that $P$ has weight 32 and is the only Weierstrass point of $X$.
All other points
	 of $X$ have gap sequence 1,2,5.  Thus $X$ is a nonclassical curve. (A curve
of
	 genus $g$ is called classical if the sequence of  Weierstrass gaps at all but
	 finitely many points of the curve is $1,2,\dots,g$.)

	 F.K. Schmidt [17] observed that there are no nonclassical smooth curves of
genus 3
	 in characteristic 2, so this example shows that singular curves can exhibit
behavior not
	 found on smooth curves.  Note that the one-parameter family of curves
	 $$\{{\cal X}_u\}= y^3z+x^4+uxz^3$$ over $k$
	 has the property that ${\cal X}_u$ is nonsingular, and hence classical, if
$u\ne 0$, but
	 ${\cal X}_0$ is nonclassical.

	 In characteristic zero, C. Widland [unpublished] has shown that every
Gorenstein curve
	 (of arithmetic genus at least two) must have at least two Weierstrass points.
	  In particular, there cannot exist a single
	 singular point that ``uses up'' all the Weierstrass weight as does the
singularity in this
	 example.  This should not be too surprising though, since there also exist
smooth curves
	 in positive characteristic that have a single Weierstrass point (cf. [6]).
	 \medskip
	 We can generalize the situation in the above example to obtain the following
result, which
	 describes the Weierstrass points on a rational curve with one unibranch
singularity
	 whose local ring is of a certain type.
         First,  we recall some facts about numerical semigroups.  A numerical
semigroup is a
	  subsemigroup of the nonnegative integers ${\bf N}$ that includes all but
	  finitely many positive integers.  The missing positive integers are called
gaps.
	    If $S$ is a numerical semigroup, then the conductor of $S$ is the least
	    integer $c$ such that  $c+{\bf N}\subseteq S$.  A numerical semigroup
	  $S$ with conductor $c$ is called symmetric
	   if  $m\in S$ if and only if $c-1-m\notin S$ for all integers $m$.  If $S$
is symmetric and
	   has $g$ gaps, then the conductor of $S$ is $2g$. E. Kunz
	   [9] showed that a unibranch curve  singularity is Gorenstein if and only if
	    the semigroup of  values associated to the singularity is a symmetric
semigroup.
	 \medskip
	   \no{\bf (2.2) Definition.}  Suppose $X$ is a rational curve and $P$ is a
	  unibranch singularity on $X$ with semigroup of values $S$.
	  Suppose the conductor of $S$ is $c$.
	  	   Let $0,n_1,n_2,\dots,n_r$ denote the nonnegative integers
	    less than $c$ that are in $S$.	   We will call $P$ a monomial unibranch
singularity
	     if the local ring $\O$ at $P$ is of the form
	   $$\O=k+kt^{n_1}+kt^{n_2} +\cdots +kt^{n_r} +t^c\tilde \O,$$
	   where $\tilde \O$ denotes the normalization of $\O$ and where $t$ is a
uniformizing
	   parameter of $\tilde\O$ that generates the function field of $X$.

	 \proclaim {(2.3) Theorem}.  Suppose $S$ is a symmetric numerical semigroup
with gaps
	 $l_1=1,l_2,\dots,$ $l_g$.  Let  $n_0=0,n_1,\dots,n_{g-1}$ be the nonnegative
	 integers less than $2g$ that are in $S$.  Suppose $X$ is a rational curve
with
	  a unique  singular point $P$ such that $P$ is a unibranch singularity and
	 the semigroup of values associated to $P$ is $S$.
	 Denote by
$\epsilon_0=0,\epsilon_1=1,
	 \epsilon_2,\dots,\epsilon_{g-1}$ the orders of $\omega_X$. Then we have:

	 \no{\it 1) $W(P)\ge \sum_{i=0}^{g-1} (n_i-\epsilon_i) +
2g\sum_{i=0}^{g-1}\epsilon_i$,
	 and equality occurs if and only if $det ({n_j\choose \epsilon_i})\not\equiv
0$
	 (mod $p$).
	 \medskip
	 \no 2) Suppose moreover that $P$ is a monomial unibranch singularity with
	 local ring
	 $$\O=k+kt^{n_1}+\cdots +kt^{n_{g-1}}+t^{2g}\tilde\O.$$  Then
	 $$\eqalign{W(P)&=\sum_{i=0}^{g-1} (n_i-\epsilon_i) +
2g\sum_{i=0}^{g-1}\epsilon_i,\cr
	 W(P_\infty)&=\sum_{i=0}^{g-1} (l_{i+1}-1-\epsilon_i),\cr}$$
	 where $P_\infty$ represents the pole of the function $t$, and there are no
other
	 Weierstrass points.}
	\medskip
	 \no {\bf Proof}. (1) Put $c=2g$.  Suppose $\sigma\in H^0(X,\omega)$.  Locally
at $P$, write
	 $\sigma=dt/h$, where $t$ is a local coordinate on the normalization of $X$
centered
	 at the point lying over $P$.
	 Suppose ord$_Ph=m+1$.  If $m\in S$, then there would be a
	  function $f\in \O_P$ such  that ord$_Pf=m$; but then  the residue of
$f\sigma$ would
	  not be zero at 0.  Therefore, $m$ must be a gap of $S$.  It follows that
there exist
	  linearly independent dualizing differentials $\tau_1,\tau_2,\dots,\tau_g$
such that
	  if we write $\tau_j=dt/h_j$, then ord$_Ph_j=l_j+1$.  Note that $l_g+1=c$.
For a
	  generator of $\omega_P$, we may take $\tau=\tau_g=dt/h_g$.  With notation as
in \S 1,
	  the functions $f_1,f_2,\dots,f_g$ used to construct the wronskian are, after
renumbering,
	  $$f_1=1,f_2= h_g/h_{g-1},\dots,f_g=h_g/h_1.$$
	  Note that we have   $$\hbox{ord}_P  f_j=c-(l_{g-j+1}+1)=n_{j-1}
	  \quad\hbox{for }j=1,2,\dots,g,$$
	  since $S$ is symmetric.  At $P$, the wronskian  vanishes to order
	  $$\eqalign{&c(\epsilon_0+\epsilon_1+\cdots +\epsilon_{g-1}) +
	   \hbox{ ord}_P\det(D^{(\epsilon_i)}f_j)\ge\cr
	   &2g\sum_{i=0}^{g-1}\epsilon_i +
\sum_{i=0}^{g-1}(n_i-\epsilon_i).\cr}\eqno{(\hbox{$*$})}$$
	   The assertion about when the equality holds follows from [19, Theorem 1.5].

	   \no (2) From the assumption that $P$ is a monomial unibranch singularity,
it follows
	   that $dt/t^{l_j+1}$ is in fact a dualizing differential, so we may take
$\tau_j=dt/t^{l_j+1}$
	   and $\tau=dt/t^c$ generates $\omega_P$.
	   At $P_\infty$, a local coordinate is $u=1/t$ and $\tau_1=-du$ generates
	   $\omega_{P_\infty}$.
	    Therefore, at $P_{\infty}$, our basis of dualizing differentials can be
written
	   $\tau_j=u^{l_j-1}\tau_1$ for $j=1,2,\dots,g$.
	    Hence the Weierstrass gap sequence  at $P_\infty$ is $l_1,l_2,\dots,l_g$.
	       Thus, we  have
	     $$W(P_\infty)\ge \sum_{i=0}^{g-1}(l_{i+1}-1-\epsilon_i ).\eqno{(\hbox{$*
*$})}$$
	     Now, if we add the right sides of equations ($*$) and ($* *$) and use the
fact
	     that $n_i+l_{g-i}=2g-1$ for each $i=0,1,\dots,g-1$, then we obtain

$$\eqalign{&2g\sum_{i=0}^{g-1}\epsilon_i+\sum_{i=0}^{g-1}(n_i-\epsilon_i)+
	     \sum_{i=0}^{g-1}(l_{i+1}-1-\epsilon_i)\cr
	     &=\sum_{i=0}^{g-1}(n_i+l_{g-i})-g+(2g-2)\sum_{i=0}^{g-1}\epsilon_i \cr
	     &=g(2g-1)-g +(2g-2)\sum_{i=0}^{g-1}\epsilon_i\cr
	     &= (\sum_{i=1}^{g-1}\epsilon_i +g)(2g-2).\cr}$$
	     Since this is the total Weierstrass weight of all Weierstrass points of
$\omega$, we
	     must have equality in ($*$) and in ($* *$) and the Theorem follows.\qed
	     \medskip
	 Notice that, for $X$ as in  part (2) of Theorem (2.3), the symmetric
semigroup
	  $S$ is the Weierstrass semigroup of nongaps
	 at the point $P_\infty$.  Thus, as was previously noted (in characteristic
zero)
	 by K.-O. St\"ohr
	 [20], it is easy to see that every numerical symmetric semigroup occurs
	 as the Weierstrass semigroup at a point on some Gorenstein curve.  It is
still unknown
	 if every such semigroup occurs as the Weierstrass semigroup at a point on
some
	 nonsingular curve.

	 \proclaim{(2.4) Corollary}.  Suppose that $X$ is as in part (2) of Theorem
(2.3).
	 If the characteristic of $k$ does not divide the integer
	 $$\prod_{i>j} (l_i-l_j)/(i-j),$$
	 then $X$ is classical.

	 \no{\bf Proof.}  The $\omega$-orders at $P_\infty$ are
	 $l_1-1,l_2-1,\dots,l_g-1$.  The Corollary then follows from the proof of
Corollary 1.7 in
	 [19].\qed

	 \medskip
	 In order to completely characterize the rational curves with one monomial
	 unibranch singularity
	 that uses up all the Weierstrass weight, we need the following result.  For
the remainder
	 of this section we assume that $k$ has characteristic $p>0$.

	 \proclaim {(2.5) Proposition}.
	   Suppose $V$ is a linear system on $X$
	 of (affine) dimension $s$.  Let $\epsilon_0,\epsilon_1,\dots, \epsilon_{s-1}$
be the order
	 sequence of $V$.  Suppose $Q\in X$ is a smooth point and let
	 $$0=\epsilon_0(Q)<\epsilon_1(Q)<\cdots <\epsilon_{s-1}(Q)$$
	 be the $(V,Q)$-orders.  Then the following are equivalent.

	 {\it (i) The $V$-Weierstrass weight of $Q$ is
$\sum_{i=0}^{s-1}(\epsilon_i(Q)-\epsilon_i)$.
	 \medskip
	 (ii) The sequence $\epsilon_0,\epsilon_1,\dots,\epsilon_{s-1}$ is the minimal
sequence,
	  in the lexicographic order, such that
	  $$\det \left( {\epsilon_j(Q)\choose \epsilon_i}\right) \not\equiv 0
\quad(\hbox{mod }p).$$
	  \medskip
	  (iii) The orders of the morphism from $\P^1$ to  $\P^{s-1}$ defined by
	  $$t\mapsto (1: t^{\epsilon_1(Q)} : t^{\epsilon_2(Q)}:\cdots :
t^{\epsilon_{s-1}(Q)})$$
	  are $\epsilon_0,\epsilon_1,\dots,\epsilon_{s-1}$.
	  \medskip
	  (iv) The sequence $\epsilon_0,\epsilon_1,\dots,\epsilon_{s-1}$ is the
minimal
	  sequence, in the lexicographic order, such that
	  $$\det \left({\epsilon_{s-1}(Q)-\epsilon_{s-1-j}(Q)\choose
\epsilon_i}\right)\not\equiv 0\quad
	  (\hbox{mod }p).$$}
	  \medskip
	  \no {\bf Proof.} The equivalence of (i) and (ii) follows from [19, Theorem
1.5 and
	  Proposition 1.6].  The equivalence of (ii) and (iii) follows from [19,
Proposition
	  1.6 and subsequent Remark].   To see  that (iv) is equivalent to the other
	  statements, first note that if $t\ne 0$, then
	  $$\eqalign{(1: t^{\epsilon_1(Q)} : t^{\epsilon_2(Q)}:\cdots :
t^{\epsilon_{s-1}(Q)})&=
	  (1: (1/t)^{-\epsilon_1(Q)} : (1/t)^{-\epsilon_2(Q)}:\cdots :
(1/t)^{-\epsilon_{s-1}(Q)})\cr
	  &=((1/t)^{\epsilon_{s-1}(Q)}:(1/t)^{\epsilon_{s-1}(Q)-\epsilon_1(Q)}:\cdots:
	    1).\cr}$$
	  	  Thus, the orders of the morphism in (iii) are the same as the orders of
the morphism
	  from $\P^1$ to $\P^{s-1}$ defined by
	  $$1/t\mapsto (1:(1/t)^{\epsilon_{s-1}(Q)-\epsilon_{s-2}(Q)}:\cdots:
(1/t)^{\epsilon_{s-1}(Q)-
	  \epsilon_1(Q)}:(1/t)^{\epsilon_{s-1}(Q)}).$$
	    The equivalence of (iv) with the other statements now follows from the
equivalence
	   of (ii) and (iii).\qed
	   \medskip
	   As one consequence of Proposition (2.5), we obtain the following general
corollary,
	   which is a new result in the theory of Weierstrass points in positive
	   characteristic.
	  \proclaim {(2.6) Corollary}. Suppose $Q$ is a smooth $V$-Weierstrass point
of $X$ with
	  the $(V,Q)$-orders being $\epsilon_i(Q), i=0,1,\dots,s-1$.  If the sequence
$\epsilon_0(Q),
	  \epsilon_1(Q),\dots,\epsilon_{s-1}(Q)$ satisfies the $p$-adic criterion,
then
	  $$W_V(Q)>\sum_{i=0}^{s-1}(\epsilon_i(Q)-\epsilon_i),$$
	  where $\epsilon_i, i=0,1,\dots,s-1$ are the orders of $V$.

	  \no{\bf Proof}.  It follows from Propositions (1.3) and (2.5) that if
$$W_V(Q)=
	  \sum_{i=0}^{s-1}(\epsilon_i(Q)-\epsilon_i),$$ then the
	  $\epsilon_i(Q)$ would be the orders of $V$.  But then $Q$ would not be a
	  $V$-Weierstrass point.
	  \line{\hfil \qed}
	  \medskip
	  \proclaim {(2.7) Corollary}.  Let $X$ be a curve as in part (2) of Theorem
(2.3).  Then the
	  $\omega$-orders $\epsilon_0,\epsilon_1,\dots,$ $\epsilon_{g-1}$ are minimal,
in the
	  lexicographic order, such that
	  $$\det\left({l_{j+1}-1\choose \epsilon_i}\right)\not\equiv 0\quad(\hbox{mod
}p).$$

	  \no{\bf Proof.} In the proof of Theorem (2.3), it was shown that
	  $\epsilon_i(P_\infty)=l_{i+1}-1$ for $i=0,1,\dots,g-1$.
	  The Corollary then follows from Theorem (2.3) and
	  Proposition (2.5), taking $V=H^0(X,\omega)$ and $Q=P_\infty$.\qed
	  \medskip
	  Finally, we obtain a characterization of rational curves with a single
	  monomial unibranch singularity that uses up all the Weierstrass weight.
	  \proclaim {(2.8) Corollary}.  Let $X$ be a curve as in  part (2) of
	   Theorem (2.3).  Then the singularity
	  $P$ uses up all the Weierstrass weight if and only if the sequence
	  $$l_1-1,l_2-1,\dots,l_g-1$$
	  satisfies the $p$-adic criterion.

	  \no{\bf Proof.}  By Theorem (2.3), $P$ uses up all the Weierstrass weight if
and
	   only if $W(P_\infty)=0$; i.e.,   if and only if
	   $l_{i+1}-1$ are the orders of $\omega$.  If these are the orders of
$\omega$, then,
	   by Proposition (1.1),
	   they satisfy the $p$-adic criterion.  Conversely, if the sequence
	   $l_1-1,l_2-1,\dots,l_g-1$ satisfies
	   the $p$-adic criterion, then by Proposition (1.3) these integers are the
	    orders of the morphism from $\P^1$ to $\P^{g-1}$ defined by
	    $$t\mapsto (1:t^{l_2-1}:\dots:t^{l_g-1}).$$
	    Therefore, from Proposition (2.5), these integers must be the orders of
$\omega$
	    and the weight of $P_\infty$ is 0.\qed
	  \medskip
	  Three other examples of curves as in part (2) of Theorem (2.3) such that
	   the singularity uses up
	  all the Weierstrass weight are:

	  (1) $p=2, g=7, S=\langle 4,6, 11\rangle=4\N + 6\N + 11\N$.

	  (2) $p=3, g=4, S=\langle 3,5\rangle= 3\N + 5\N$.

	  (3) $p=5, g=6, S=\langle 4,5\rangle=4\N + 5\N$.
	  \medskip
	  We conclude this section with an example of a rational curve with a single
	  (nonmonomial) unibranch  singularity and two smooth Weierstrass points.
	  \medskip
	  \no {\bf (2.9) Example.}  Let $X$ be the rational Gorenstein curve obtained
from
	  the projective line over $k$ by replacing the local ring $\tilde\O$ at 0
with the
	  ring
	  $$\O=k+k(t^3+t^5) +kt^4 + t^6\tilde\O,$$
	  where $t$ is a uniformizing parameter of $\tilde\O$ that generates the
	  function field of $X$. Let $P$ denote the singular point of
	  $X$.  The semigroup of values at $P$ is $3\N +4\N$, the same as the
	  semigroup of values in Example (2.1).  But notice here that $dt/t^6$ is not
	 a dualizing differential since
	 $$\hbox{Res}_0(t^3+t^5) \> {dt\over t^6}=1.$$
	 It is not hard to see that a basis for $H^0(X,\omega)$ is
	 $$\tau_1=(1-t^2) dt/t^6, \tau_2=dt/t^2, \tau_3=dt/t^3$$
	 and $\tau_1$ generates $\omega_P$. One may check easily that the order
	  sequence at infinity is 0,1,2, so $X$ is classical (in all characteristics).

	 If the characteristic is not 2, then the wronskian, on the open subset of
	 $X$ obtained by excluding the point at infinity and the zeros of $1-t^2$, is
	 $${t^{22}(6-t^2)\over (1-t^2)^6},$$ up to a nonzero constant.
	 Thus, if the
characteristic is not 2 or 3,
	  $P$ has weight 22 and there are two weight one
	 Weierstrass points at the two square roots of 6.  If the characteristic is 3,
then
	 $P$ has weight 24 and is the unique Weierstrass point of $X$.
	  A computation of the wronskian in characteristic 2 shows that $P$  has
	  weight 24 in this case as well.

	 \goodbreak
	 \medskip
	 \no {\bf 3.} In this section, we will consider rational curves with several
	  unibranch singularities.
	 We will assume for the remainder of this article that $k$ has characteristic
0.
	   We do this so that our formulas are not
	 overly complicated.  One can modify these formulas, taking into account the
order
	 sequences of the canonical bundles involved, and obtain analogous results
	 in positive characteristic.

	 \medskip
	     \no {\bf (3.1) Definition}.  If $S$ is a numerical semigroup with gaps
$l_1,l_2,\dots,
	     l_\delta$, then we define the {\it weight of }$S$, denoted $wt(S)$, by
	     $$wt(S) =\sum_{j=1}^\delta (l_j-j).$$
	     We define the weight of ${\bf N}$ to be 0.
	     \medskip
	     \no {\bf (3.2) Remark}.  We always have $l_\delta\le 2\delta-1$ (and the
equality occurs if
	     and only if $S$ is symmetric).  To see this, note that if $c$ is the
conductor of $S$
	     and if $x+y=c-1=l_\delta$, then either $x$ or $y$ must be a gap.  Hence
among the
	     nonnegative integers less than $c$, there are at least as many gaps as
there are
	     nongaps.
	     \proclaim {(3.3) Lemma}.  If $S$ is a numerical semigroup with $\delta$
gaps and if
	     $0=n_0,n_1,\dots,n_{\delta-1}$ are the elements of $S$ less than
$2\delta$, then
	     $$\sum_{i=0}^{\delta-1}(n_i-i)=(\delta-1)\delta -wt(S).$$

	     \no{\bf Proof.} Let $l_1,l_2,\dots,l_\delta$ denote the gaps of $S$.  We
have
	     $$\sum_{j=1}^\delta l_j
+\sum_{i=0}^{\delta-1}n_i=\sum_{k=0}^{2\delta-1}k=
	     \delta(2\delta-1).$$
	     Hence
	     $$\eqalign{wt(S)&=\sum_{j=1}^\delta l_j-\delta(\delta+1)/2\cr
	                               &= \delta(2\delta-1)-\sum_{i=0}^{\delta-1}n_i
-\delta(\delta+1)/2\cr
				       &=
\delta(2\delta-1)-\sum_{i=0}^{\delta-1}n_i-(\delta^2-(\delta-1)\delta/2)\cr
				       &=(\delta-1)\delta -\sum_{i=0}^{\delta-1}(n_i-i).\qed\cr}$$
	     \medskip
	   The next theorem treats the weight of
	 a unibranch singularity on an arbitrary (in particular, not necessarily
rational)
	 Gorenstein curve.  This theorem generalizes a result of C. Widland [21] in
	  the case of a simple cusp.

	  We recall that if $P$ is a singular point of $X$, then by the partial
normalization of $X$
	  at $P$ one means the curve obtained from $X$ by desingularizing only the
singularity
	  $P$.

	 \proclaim {(3.4) Theorem}.  Suppose $X$ is a (not necessarily rational)
Gorenstein curve
	 of arithmetic genus $g$.  Suppose $P\in X$ is a unibranch singularity.  Put
$\delta=\delta_P$.
	 Let $S$ denote the semigroup of values at $P$.  Let $Y$ denote the partial
normalization
	 of $X$ at $P$ and let $Q$ denote the point of $Y$ that lies over $P$.  Then
	 $$W_X(P)=\delta (g-1)(g+1) -wt(S) + W_Y(Q).$$

	 \no {\bf Proof.}  Let $\sigma_1,\sigma_2,
	 \dots,\sigma_{g-\delta}$ be a basis of $H^0(Y,\omega_Y)$.  Locally at $Q$,
write
	 $\sigma_i=f_i\> dt$, for $i=1,\dots,g-\delta$, where $t$ is a local
coordinate on $Y$
	 centered at $Q$.  Put $r_i=\hbox{ord}_Qf_i$.  We will assume that the basis
of differentials
	 has been chosen so that $r_1<r_2<\cdots<r_{g-\delta}$.

	 Dualizing differentials on $Y$ are also dualizing differentials on $X$, so we
may extend
	 the above differentials to  a basis
	 $$\tau_1,\tau_2,\dots,\tau_\delta,\sigma_1,\sigma_2,\dots,\sigma_{g-\delta}$$
	 of $H^0(X,\omega_X)$.  Let $l_1,l_2,\dots,l_\delta$ denote the gaps of the
semigroup of
	 values at $P$.
	  As in Theorem (2.3), we may assume that, locally at $P$, we have
	 $\tau_i=dt/h_i$, where ord$_Ph_i=l_i+1$ for $i=1,2,\dots,\delta$.  Since the
semigroup
	 of values at $P$ is symmetric, we have $l_\delta+1=2\delta=c$, the conductor
of the
	 semigroup of values.  The differential $\tau_\delta$ generates
$\omega_{X,P}$.  Put
	 $h=h_\delta$.

	  Since we are assuming the characteristic is 0, the orders of $\omega_X$ are
	  $0,1,\dots,g-1$ and the sum of the orders is $N=(g-1)g/2$.   As in the proof
of
	  Proposition (1.6),  we have
	   $$W_X(P)=\delta(g-1)g + \hbox{ord}_Q
W_t(1,h/h_{\delta-1},h/h_{\delta-2},\dots,h/h_1,
	   hf_1,hf_2,\dots,hf_{g-\delta}),$$
	   where $W_t$ denotes the ordinary wronskian (obtained by differentiating
with respect
	   to $t$) of the given functions and where we have used the fact that the
order of a
	   function in $\O_P$ is the same at $P$ as it is at $Q$.  Notice that each of
the functions
	   $1,h/h_{\delta-1},\dots,h/h_1,$ $hf_1,\dots,hf_{g-\delta}$ has a different
order at $Q$.
	   Indeed, we have ord$_Qh/h_{\delta-i}=c-(l_{\delta-i}+1)=n_i$,
	   where $n_0,n_1,\dots,n_{\delta-1}$ are the elements in $S$ that are less
than $2\delta$,
	    and ord$_Qhf_j=c+r_j=2\delta+r_j$.
	     Hence the order of the wronskian of these functions at $Q$ may be easily
computed
	     as in [3, p.82].  We have
	     $$\eqalign{\hbox{ord}_QW_t(1,h/h_{\delta-1},\dots,hf_{g-\delta})&=
	     \sum_{i=0}^{\delta-1}(n_i-i) + \sum_{j=1}^{g-\delta}
(2\delta+r_j-(\delta+j-1))\cr
	     &=\sum_{i=0}^{\delta-1}(n_i-i)+\delta(g-\delta)
+\sum_{j=1}^{g-\delta}r_j-j+1\cr
	     &=\sum_{i=0}^{\delta-1}(n_i-i)+\delta(g-\delta) + W_Y(Q).\cr}$$
	     The Theorem now follows from Lemma (3.3).\qed

	     \medskip
	     Using Theorems (2.3) and (3.4), one may compute the weight of each
singularity on a
	     rational curve with precisely two monomial unibranch singularities.

	     \proclaim {(3.5) Theorem}.  Suppose $X$ is a rational Gorenstein curve
	     of arithmetic genus $g$ with two monomial unibranch
	     singularities $P_1$ and $P_2$. For $i=1,2$ put $\delta_i=\delta_{P_i}$,
let $S_i$
	      denote the semigroup of values at $P_i$, and let $t_i$ denote the
uniformizing
	      parameter at $P_i$ that generates the function field of $X$ used to
define the
	      local ring at $P_i$ as in Definition (2.2).
	      Let $Q_1$ and $Q_2$ be the points on ${\bf P}^1$ lying over $P_1$ and
$P_2$,
	      respectively.  Then

	     \no{\it 1)  If $Q_1$ and $Q_2$ are the poles  of  the functions $t_2$ and
	     $t_1$, respectively (i.e., if $t_1t_2$ is a nonzero constant), then
	      $$\eqalign{W_X(P_1)&=\delta_1(g-1)(g+1)
	      -wt(S_1)+wt(S_2)\cr
	       W_X(P_2)&= \delta_2(g-1)(g+1) -wt(S_2)+
	       wt(S_1),\cr}$$ and there
	     are no smooth Weierstrass points on $X$.

	      \no 2) If $Q_1$ is the pole of the function $t_2$, but $Q_2$ is not
	      the pole of the function $t_1$, then
	      $$\eqalign{W_X(P_1)&=\delta_1(g-1)(g+1)-wt(S_1)+wt(S_2)\cr
	      W_X(P_2)&=\delta_2(g-1)(g+1)-wt(S_2),\cr}$$
	      and there are $wt(S_1)$ smooth Weierstrass points on $X$, counting
	      multiplicities.

	      \no 3) If $Q_1$ and $Q_2$ are not the poles of the functions $t_2$ and
	      $t_1$, respectively, then
	      $$\eqalign{W_X(P_1)&=\delta_1(g-1)(g+1)
	      -wt(S_1)\cr
	      W_X(P_2)&= \delta_2(g-1)(g+1) -wt(S_2),\cr}$$
	       and there are $wt(S_1)+wt(S_2)$ smooth Weierstrass points on $X$,
	        counting multiplicities.}
		\medskip
		\no{\bf Proof.} By Theorem (3.4)  we have
		 $$W_X(P_1)=\delta_1(g-1)(g+1)-wt(S_1)+
		 W_{Y_1}(Q_1),$$
		 where $Y_1$ is the partial normalization of $X$ at $P_1$ and $Q_1$ is the
point on
		 $Y_1$ that lies over $P_1$.  Now, $Y_1$ is a rational curve with the unique
	 monomial unibranch singularity
		 $P_2$.  Hence we see from Theorem (2.3) that $W_{Y_1}(Q_1)=wt(S_2)$
		  if $Q_1$ is the pole of the function $t_2$
		  and is 0 otherwise.  A similar argument holds with regard to
		 $W_X(P_2).$  The assertions about the number of smooth Weierstrass
		  points on $X$
		 follow by adding the weights of $P_1$ and $P_2$ and subtracting from
		 $g^3-g$, which is the total of all the weights.  Note that
$g=\delta_1+\delta_2$.\qed
		 \medskip
		 \no{\bf  (3.6) Example}. Suppose that $X$ is the rational curve obtained
		  from ${\bf P}^1$
		 by creating two monomial unibranch singularities $P_0$ and $P_1$, each
with
		 semigroup of values generated by 3 and 4,  at 0 and 1.  Then $X$ has
		 arithmetic genus 6 and the Weierstrass weight of each singularity is 103.
The
		 total Weierstrass weight is 210, and it may be seen, by computing the
wronskian on the
		 smooth locus of $X$, that there are four distinct smooth Weierstrass points
(each of
		 weight one).  We note that the point at infinity  is not a Weierstrass point
on $X$, but it
		 is a Weierstrass point on $Y_1$, the partial normalization of $X$ at $P_1$.
 		 The existence of a function with a zero of order 3 at $P_0$ shows that 3 is
a nongap at
		 infinity on the curve $Y_1$, but this function is not regular at $P_1$ on
$X$ and 3 is
		 not a nongap at infinity on $X$.
		 \medskip
		 The situation in part (3) of Theorem  (3.5) may be generalized as follows.
		 \proclaim {(3.7) Theorem}.  Suppose that $X$ is a rational Gorenstein
curve
		 of arithmetic genus $g$ with
		 unibranch singularities $P_1,P_2,\dots,P_n$ as its only singularities.
		 For $i=1,2,\dots,n$, let $S_i$  denote the semigroup of values at $P_i$,
		  and assume that the
		 point $Q_i$ that lies over $P_i$ is not a Weierstrass point of
		 the partial normalization $Y_i$ of $X$ at $P_i$.
		Then the number of smooth Weierstrass points on $X$, counting
multiplicities,
is
		$\sum_{i=1}^n wt(S_i)$.

		\no{\bf Proof.}  Since by hypothesis $W_{Y_i}(Q_i)=0$, it follows from
Theorem (3.4)
		that the number of smooth Weierstrass points on $X$, counting
		multiplicities, is
		$$\eqalign{g^3-g-\sum_{i=1}^nW_X(P_i)&=
		g^3-g-\sum_{i=1}^n (\delta_i(g-1)(g+1)-wt(S_i))\cr
				&= \sum_{i=0}^n wt(S_i),\cr}$$
                 since $g=\sum_{i=0}^n\delta_i$.\qed
		 \medskip
		 \goodbreak
		 \no {\bf 4.} We now consider singularities with two branches.  Suppose
$X$
is an
		 integral, projective Gorenstein curve over an algebraically closed field of
		 characteristic zero.
		 Suppose $P\in X$ is a singularity with two branches.  Let $\pi:Y\rightarrow
X$
		 denote the partial normalization of $X$ at $P$
		  and let $Q_1,Q_2\in Y$ denote the
		 points lying over $P$.  Let $\nu_1$ and $\nu_2$ denote the discrete
		  valuations associated to $Q_1$ and $Q_2$, respectively.  The value
semigroup $S$
		  at  $P$ is given by
		  $$S=\{(\nu_1(f),\nu_2(f))\in {\bf N}\times{\bf N} : f\in\O_P,f\ne 0\}.$$
		  Such semigroups have been studied by, among others, the first author [4]
		   in the case of plane  curves and by F. Delgado
		   in the cases of plane curves [1] and Gorenstein curves [2]
		  with an arbitrary number of branches.  Let $\xi=(\xi_1,\xi_2)$ denote the
conductor of
		  $S$; i.e., $\xi$ is the minimum element in $S$, with respect to the
		  product order on $\N\times \N$, such that  $\xi+ \N\times\N \subseteq S$.

		   We think of the semigroup $S$ as being a set of points in the plane.
		   Let $S_i=\pi_i(S)$, for $i=1,2$, denote the projections of $S$ onto the
coordinate axes.
		    For $i=1,2$,  let $\delta_i$ denote the number of gaps of
		   $S_i$
		   and let  $c_i$ denote the conductor of $S_i$.  We note that $S_1$ and
$S_2$ need
		   not be symmetric semigroups (see Example (4.10) below).

		    Since $\O=\O_P$ is a Gorenstein ring, the semigroup $S$ also has
		     certain symmetry properties, which  we now recall.
		   \medskip
		   \no {\bf (4.1) Definitions.}  For $x\in \N$, the vertical fiber at $x$,
denoted VF$(x)$, is
		   defined by $$\hbox{VF}(x)=\{(x,y')\in \N\times\N : (x,y')\in S\}.$$
		   For $y\in \N$, the horizontal fiber at $y$, denoted HF$(y)$, is defined
similarly.
		   A point $(x_1,y_1)$ is said to be above (resp. to the right of ) another
point
		   $(x_2,y_2)$ if $x_1=x_2$ and $y_1>y_2$ (resp. if $y_1=y_2$
		    and $x_1>x_2$).
		   Put
		   $$\Delta((x,y))=\{(x',y')\in S: (x',y') \hbox{ is either above or to the
right of }
		   (x,y)\}.$$
		   A point $(x,y)\in S$ is called a {\it maximal point}  (or simply a {\it
maximal})
		   if $\Delta((x,y))=\emptyset$.
		   \proclaim {(4.2) Lemma}.

		   {\it\no 1) Suppose $n\in S_1$ and $n<\xi_1$.  Then
		   $$(n,\xi_2)\in S\Leftrightarrow \xi_1-1-n\notin S_1\Leftrightarrow
\hbox{VF($n$) is
		   infinite.}$$
		   \no 2) Suppose $n\in S_2$ and $n<\xi_2$.  Then
		   $$(\xi_1,n)\in S\Leftrightarrow \xi_2-1-n\notin S_2\Leftrightarrow
		   \hbox{HF($n$) is infinite.}$$}
		   \medskip
		   \no{\bf Proof}. [2, Lemma (1.8) and Theorem (2.3)]\qed
		   \medskip

		   Put $\mu=(\xi_1-1,\xi_2-1)$.   From [2, Corollary (2.7)], we have that
		   $\mu$ is a maximal point in $S$.  This point plays a role in $S$ analogous
to the
		   number $c-1$ in a symmetric numerical semigroup.
		   More precisely,  one has the following result.

		   \proclaim {(4.3) Proposition}. (Symmetry properties of $S$).  The
		   semigroup $S$ has the following symmetry properties.

		   {\it \no1) For any $(x,y)\in \Z\times\Z$,
		   $$(x,y)\in S\Leftrightarrow \Delta(\mu-(x,y))=\emptyset.$$
		   \medskip
		   \no 2) For any $(x,y)\in \N\times\N$, $$(x,y)\hbox{ is a maximal of
}S\Leftrightarrow
		   \mu-(x,y)\hbox{ is a maximal of } S.$$}
		   \medskip
		   \no{\bf Proof}.  Delgado [2, Theorem (2.8)] establishes the first property
above
		   and notes ([2, Remark (2.9)]) that the second property also holds.\qed

		   \proclaim {(4.4) Lemma}.

		   \no{\it 1) Suppose $n<\xi_1$.  Then VF$(n)$ is infinite if and only if
		   $\xi_1-1-n$ is a gap of $S_1$.

		   \no 2) Suppose $n<\xi_2$.  Then HF$(n)$ is infinite if and only if
$\xi_2-1-n$ is a
		   gap of $S_2$.}
		   \medskip
		   \no {\bf Proof.}
		     Suppose VF$(n)$ is infinite, with $n<\xi_1$.
		       Then there exists a point $(n,y)\in S$ with
		     $y> \xi_2$.  By adding the function corresponding to this point with the
		     function corresponding to $(\xi_1,\xi_2)$, we see that $(n,\xi_2)\in S$.
		     Therefore, by  Lemma (4.2), we have that $\xi_1-1-n\notin S_1$.
Conversely,
		     if $n<\xi_1$ and $\xi_1-1-n\notin S_1$, then we claim that $n\in S_1$.
For consider
		     the point $\alpha=(n,\xi_2-c_2)$.  Then $\mu-\alpha=(\xi_1-1-n,c_2-1)$.
		       But $\xi_1-1-n\notin S_1$ and $c_2-1\notin S_2$, so
$\Delta(\mu-\alpha)=\emptyset$.
		       It follows from Proposition (4.3) that $\alpha\in S$ and so $n\in
S_1$.
		       Hence, if $n<\xi_1$ and $\xi_1-1-n\notin S_1$, we can conclude from
Lemma
		       (4.2) that VF$(n)$ is infinite.

		        The proof of (2) is similar.\qed
		   \proclaim {(4.5) Proposition}.  The symmetry properties in Proposition
(4.3)
		   are equivalent.

		   \no {\bf Proof.} $(1)\Rightarrow (2)$: By the symmetrical form
		    of statement (2) in Proposition (4.3), it suffices to show that
		   if $(x,y)$ is a maximal of $S$, then $\mu-(x,y)$ is also a maximal of $S$.
 Now,
		   if $(x,y)$ is a maximal of $S$, then $(x,y)\in S$ and
$\Delta((x,y))=\emptyset $.  But
		   then, by applying (1) of Proposition (4.3) in both directions,
		    we see that $\Delta(\mu-(x,y))=\emptyset$ and
		   $\mu-(x,y)\in S$.  Therefore, $\mu-(x,y)$ is a maximal of $S$.

		   $(2)\Rightarrow (1)$: Assume $(x,y)\in S$.  Suppose there exists a point
in $S$
		   above $\mu-(x,y)$.  (A similar argument applies if there exists a point in
$S$ to
		   the right of $\mu-(x,y)$.) But then $(x,y)+(\mu_1-x, \mu_2-y+z)$, for some
$z>0$,
		   is in $S$, contradicting the fact that $\mu$ is a maximal point.

		   Conversely, suppose that $\Delta(\mu-(x,y))=\emptyset$.  Since
VF($\mu_1-x$)
		   and HF($\mu_2-y$) are then finite, it follows, from Lemma (4.4), that
$x\in S_1$ and
		   $y\in S_2$.  We have then two possibilities: either VF$(\mu_1-x)$ is
empty
or
		   nonempty.  If VF$(\mu_1-x)$ is empty, it follows from Lemma (4.4) that
VF$(x)$
		   is infinite.  If VF$(\mu_1-x)$ is nonempty, then the maximal point of this
fiber, call
		   it $(\mu_1-x,z)$, satisfies $z\le \mu_2-y$.  Hence, applying (2) of
		   Proposition (4.3), we have a maximal point of the form $(x,\mu_2-z)$ with
		   $\mu_2-z\ge y$.  So, in any case, one has a point in the semigroup $S$ of
the
		   form $(x,y')$ with $y'\ge y$.  Similarly, one has a point in $S$ of the
form
		   $(x',y)$ with $x'\ge x$.  If $y'=y$ or $x'=x$, we are finished.  We can
then assume
		   that $y'>y$ and $x'>x$.
		     Then the sum of the functions in the local ring $\O_P$ corresponding to
		   $(x',y)$ and $(x,y')$ is a function $f$ satisfying $\nu_1(f)=x$ and
$\nu_2(f)=y$,
		   showing that $(x,y)\in S$.\qed

			\medskip
		       Put $I$ equal to the number of maximal points in $S$.
		      From [4], we have that if $X$ is a plane curve, then
		      $I$ is also equal to the intersection number of the two branches and
the conductor
		      of $S$ is $(I+2\delta_1,I+2\delta_2)$.  We now show that these results
also hold for
		      any Gorenstein curve.

		      \proclaim{(4.6) Proposition}.  The coordinates of the conductor of $S$
are
		       $\xi_1=I+2\delta_1, \xi_2=I+2\delta_2.$

		     \no{\bf Proof.}  Consider the vertical fibers VF$(x)$ for $0\le x<
\xi_1$.
		     We will count how many of these fibers are infinite, empty, or finite
and
		     nonempty.  From Lemma (4.4), we see that the number of these vertical
		     fibers that are infinite is $\delta_1$.
		     The number of empty vertical fibers is also
		     equal to $\delta_1$.  The number of nonempty finite fibers is equal to
$I$, the number
		     of maximal points.
		    		     Therefore, $\xi_1=I+2\delta_1$.  A similar argument
		     using horizontal fibers shows that $\xi_2=I+2\delta_2$.\qed

		  \proclaim{(4.7) Corollary}. $\delta=\delta_P=I+\delta_1+\delta_2$ and $I$
is
		  the intersection number of the two branches at $P$.

		  \no {\bf Proof.} From Proposition (4.6) and the fact that $\O_P$ is
		  Gorenstein, we have $2\delta_P=2I+2\delta_1+2\delta_2$.  Therefore,
$\delta_P=
		  I+\delta_1+\delta_2$.  It then follows from
		  [7, Proposition 4] that $I$ is the intersection number of the two branches
at $P$.
		  \qed

		     \proclaim{(4.8) Corollary}.  Suppose the maximal points of $S$ are
		       $$(a_0,b_0),(a_1,b_1),\dots,(a_{I-1},b_{I-1}).$$
		       Then we have
		       $$\eqalign{\sum_{i=0}^{I-1}a_i&=I(I-1)/2 + \delta_1I\cr
		       \sum_{i=0}^{I-1}b_i&=I(I-1)/2 + \delta_2I.\cr}$$

		       \no{\bf Proof}. By Propositions (4.3) and  (4.6), $a_i$ is the first
		        coordinate of a maximal point if
		       and only if $I+2\delta_1-1-a_i$ is also the first coordinate of a
		        maximal.  Hence we have
		       $$\sum_{i=0}^{I-1}a_i=I(I+2\delta_1-1) - \sum_{i=0}^{I-1}a_i,$$
		       and the first equality in the statement of the Corollary follows.  A
similar argument
		       applies to the second coordinates of the maximal points.\qed
		 \medskip
		    Consider the rectangle  (with one vertex deleted)
		   $$R=\{(x,y)\in \N\times\N : x\le \xi_1\hbox{ and }y<\xi_2\hbox{ or }
		   x<\xi_1\hbox{ and }y\le \xi_2\}.$$
		     It follows from Lemma (4.4) that the points in $S$ that are on
		     the top edge of
		   $R$ are of the form $(\xi_1-1-l,\xi_2)$, where $l$ is a gap of $S_1$, and
the points in
		   $S$ that are on the right edge of $R$ are of the form
$(\xi_1,\xi_2-1-l')$, where $l'$ is
		   a gap of $S_2$.  If $S_1$ and $S_2$ are symmetric, then one can write
these
		   points in a nicer form.

		 \proclaim{(4.9) Corollary}.

		 \no {\it 1) If $S_1$ is symmetric and if  $m_0,m_1,\dots
		   m_{\delta_1-1}$ are the elements in $S_1$ that are less than $c_1$,
		    then the points that are in $S$ and on the upper edge of the rectangle
$R$ are
		    the points $(I+m_j,\xi_2), j=0,1,\dots,\delta_1-1$.
		  \medskip
		   \no  2)  If $S_2$ is
		    symmetric and if $n_0,n_1,\dots,n_{\delta_2-1}$ are the elements
		     in $S_2$ that are less than  $c_2$, then the points that are in $S$
		     and on the right edge of $R$ are the points $(\xi_1,I+n_k), k=0,1,
		     \dots,\delta_2-1$.}

		     \no {\bf Proof}.
		     Suppose $S_1$ is symmetric.  Then $c_1=2\delta_1$ and
		     $c_1-1-n\in S_1$ if and only if $n\notin S_1$.  Therefore, from Lemma
(4.4)
		     and Proposition (4.6),
		     $$n\notin S_1\Leftrightarrow \hbox{VF}(I+c_1-1-n)\hbox{ is infinite}.$$
		     Thus, the infinite vertical fibers are of the form VF$(I+m)$, where
$m\in S_1$.
		    	The analogous statement for $S_2$ is proved similarly.\qed
		      \medskip
		      We thank Professor K.-O. St\"ohr for the following example of a
two-branch
		      Gorenstein singularity having one branch that is not Gorenstein.
		      \medskip
		      \no{\bf (4.10) Example}.  Let $H$ be a numerical semigroup with $g$
gaps
		      $l_1,l_2,\dots,l_g$ such that $l_g=2g-2$. Clearly, $(g-1)\notin H$.
Take
		      $S_1=(g-1)\N + H$.  We claim that
		      $$S_1=H\cup \{g-1, 2g-2\};$$
		      i.e., that $S_1$ has $g-2$ gaps.  In fact, if $\alpha$ belongs to $S_1$
but not
		      to $H\cup\{g-1,2g-2\}$, then we can write $\alpha=(g-1)+h$, for some
$h\in H$.
		      Since $\alpha\ne g-1$ and $\alpha$ is a gap of $H$ we have, from [15,
Prop.
		      1.2], that
		      $$2g-2-\alpha=g-1-h=h_1\in H.$$
		      Hence $g-1=h+h_1\in H$, a contradiction.

		      Take $\O\subseteq\tilde\O=k[[t]]\times k[[u]]$ to be the local ring
given by:
		      $$\O=\left\{\left( \sum_{i=0}^\infty a_it^i, \sum_{i=0}^\infty
b_iu^i\right) :
		      a_0=b_0, a_{g-1}=b_1,  a_{2g-2}=b_2, \hbox{ and }a_i=0 \hbox{ for all
}
		      i\notin S_1\right\}.$$
		      Since $S_1$ has $(g-2)$ gaps, we have
		      $$\delta=\dim_k \tilde\O/ \O = (g-2) +3=g+1.$$
		      Clearly, the conductor ideal $\C$ of $\O$ in $\tilde \O$ is
		      $$\C=t^{2g-1}\, k[[t]]\times u^3\, k[[u]]$$ and
		      dim $\tilde \O/\C = 2g-1+3=2g+2$.  This shows that $\O$ is Gorenstein.
		      However, the semigroup of the first branch, namely $S_1$, is not
symmetric
		      if $3\notin H$.  In fact, if $S_1$ is symmetric, then its conductor
$c_1$
		      satisfies  $c_1=2(g-2)=2g-4$.  Hence 2g-5 is a gap of $S_1$ and
		      of $H$.  Again by
		      [15, Prop. 1.2], we have $2g-2-(2g-5)=3\in H$.
		      \medskip

		     We now want to describe how to find a basis of dualizing differentials
on $X$ that
		     have certain orders at $P$.  This process is analogous to finding a
``Hermitian''
		     basis  (or basis ``adapted to a point'') of regular differentials at a
		      point on a smooth curve.
		     We want to show that we can choose linearly independent
		     dualizing differentials on $X$ whose orders at $P$ are related to the
maximal points
		     of $S$ and the points of $S$ that lie on the upper edge and right edge
of
		     the rectangle $R$.
		     These differentials will be those dualizing differentials in  a
``Hermitian'' basis at $P$
		     that are not regular on $Y$ either at $Q_1$ or at $Q_2$  (or at both
points if the
		     differential corresponds to a maximal point of $S$).

		     We will use the following Riemann-Roch Theorem for zero-dimensional
subschemes
		     on a Gorenstein curve, which was proved in [22] (also cf. [5]).
		     If $J$ is a proper ideal of $\O_P$, we let $\I(J)$ denote the sheaf of
$\O_X$-ideals
		     defined by $\I(J)_P=J$ and $\I(J)_Q=\O_Q$ for all $Q\ne P$.  Put
		     $$\eqalign{h(J)&=\dim_k \hbox{Hom}_{\O_X}(\I(J),\O_X)\cr
		     \iota(J)&=\dim_kH^0(X,\I(J)\otimes\omega)\cr
		     d(J)&=\dim_k\O_P/J.\cr}$$
		     The elements of Hom$_{\O_X}(\I(J),\O_X)$ may be identified with
rational
funtions
		     $f$ on $X$ such that $fJ\subseteq \O_P$ and $f\in \O_Q$ for all $Q\ne
P$.

		     \proclaim {(4.11) Theorem}. $h(J)-\iota(J)=d(J)+1-g.$

		     Let $C$ denote the conductor of $\O_P$ in its normalization
$\tilde\O_P$.

		     \proclaim {(4.12) Lemma}. $h(C)=1$.

		     \no{\bf Proof.} This follows from the fact that
Hom$_{\O_P}(C,\O_P)=\tilde\O_P$
		      (cf. the proof of Proposition (2.2) of [14]).\qed

		     \proclaim {(4.13) Proposition}.  Suppose $\tau\in H^0(X,\omega)$
generates
		     $\omega_P$.  Suppose
		      that
		     $$\O_P=J_0\supset J_1\supset \cdots \supset J_{\delta-1}\supset
J_{\delta}=C$$
		     is a strictly decreasing chain of $\O_P$-ideals.  Then there exist
$\delta$ linearly
		     independent dualizing differentials
$\tau_1,\tau_2,\dots,\tau_{\delta}\in H^0(X,\omega)$
		     such that, locally at $P$, we have $\tau_i=f_i\tau$ with $f_i\in
J_{i-1}\setminus
		     J_i$ for $i=1,2,\dots,\delta$.

		     \no{\bf Proof}.  Note that $d(J_i)=i$ since $\dim_k\O_P/C=\delta$.
		       Since $J_i\supseteq C$
		     and $h(C)=1$, it follows that $h(J_i)=1$ for all $i$.  Therefore, from
Theorem (4.11),
		     we have $\iota(J_i)=\iota(J_{i-1})-1$ for $i=1,2,\dots,\delta$.  Thus,
there exists
		     $\tau_i\in H^0(X,\I(J_{i-1})\otimes\omega)\setminus
H^0(X,\I(J_i)\otimes\omega)$
		     for $i=1,2,\dots,\delta$.
		       Then, locally at $P$, we have $\tau_i=f_i\tau$ for some $f_i\in
J_{i-1}\setminus
		       J_i$.\qed
		       \medskip
		       \no {\bf (4.14) Definition.}  For $(x,y)\in \N\times\N$, put
		       $$J(x,y)=\{f\in \O_P : \nu_1(f)\ge x\hbox{ and }\nu_2(f)\ge y\}.$$
		       \medskip
		       \no{\bf (4.15) Definition.} Suppose $\sigma\in H^0(X,\omega)$ and
$\tau$
		       generates $\omega_P$.  Locally at $P$, write $\sigma= f\tau$, with
$f\in \O_P$.
		         Then put $\nu_1(\sigma)=\nu_1(f)$ and $\nu_2(\sigma)=\nu_2(f)$.

		       \proclaim{(4.16) Theorem}. There
		        exist $\delta$ linearly independent dualizing
		       differentials $$\tau_0,\tau_1,\dots,\tau_{\delta-1}\in H^0(X,\omega)$$
such that

		       \no{\it 1) For each maximal point $(a,b)\in S$,
		        there is a  $\tau_i, 0\le i\le I-1,$ such that
		      $\nu_1(\tau_i)=a$, and $\nu_2(\tau_i)=b$.
		       \medskip
		       \no 2)  For each point in $S$ of the form $(r,\xi_2)$, with $r<\xi_1$,
		       there is a $\tau_j,I\le j\le I+\delta_1-1,$ such that
		      $ \nu_1(\tau_j)=r$, and $\nu_2(\tau_j)\ge \xi_2$.
		       \medskip
		       \no 3) For each point in $S$ of the form $(\xi_1,s)$, with $s<\xi_2$,
                        there is a $\tau_k, I+\delta_1\le k\le \delta-1,$ such
that
			 $\nu_1(\tau_k)\ge \xi_1$ and $\nu_2(\tau_k)=s$.}
		       \medskip
		       \no{\bf Proof.}  Let $0=x_1,x_2,\dots,x_{I+\delta_1}$ be the
nonnegative integers
		       such that $ x_k<\xi_1$ and $VF(x_k)\ne\emptyset$.
		       Let $(\xi_1,s_0),(\xi_1,s_1),\dots,(\xi_1,s_{\delta_2-1})$ denote the
points in $S$
		       on the right edge of the rectangle $R$.  (These points correspond to
infinite
		       horizontal fibers that lie below the line $Y=\xi_2$.)
		       Let $C$ denote the conductor of $\O_P$ in its normalization. Consider
the following
		       chain of ideals in $\O_P$:
		       $$\eqalign{&J(0,0)\supset J(x_2,0)\supset\cdots\supset
J(x_{I+\delta_1},0)
		       \supset\cr
		       &J(\xi_1,s_0)\supset J(\xi_1,s_1)\supset\cdots\supset J(\xi_1,
		       s_{\delta_2-1})\supset C.\cr}\eqno{( * )}$$
		       This is a proper chain of  ideals as in Proposition (4.13).
		         Hence there exist $\delta$
		       linearly independent dualizing differentials $\sigma_0,\sigma_1,\dots,
		       \sigma_{\delta-1}$
		       as in Proposition (4.13).   The last $\delta_2$ of these
differentials,
		       call them $\tau_{I+\delta_1},\dots,\tau_{\delta-1}$, satisfy condition
		       (3) in the statement of the Theorem.

		       In a similar manner, we may find $\delta_1$ differentials, call them
		       $\tau_I,\dots,\tau_{I+\delta_1-1}$, satisfying condition (2) in the
		       statement of the Theorem.

		       Suppose the maximal points of $S$ are
		       $$(a_0,b_0)<(a_1,b_1)<\cdots < (a_{I-1},b_{I-1}),$$
		       ordered lexicographically.   One of the differentials, call it
$\sigma$,
		       that  we found using the chain $( * )$ above satisfies
		       $$\nu_1(\sigma)=a_{I-1}, \nu_2(\sigma)\le b_{I-1}.$$
		       If $\nu_2(\sigma)\ne b_{I-1}$, then $\nu_2(\sigma)=s_k$ for some
		       $k, 0\le k\le \delta_2-1$.  In that case, a suitable linear
combination of
		       $\sigma$ and the differential $\tau_{I+\delta_1+k}$ will yield a
differential
		       $\bar\sigma$ such that
		       $$\nu_1(\bar\sigma)=a_{I-1}\hbox{ and
}\nu_2(\bar\sigma)>\nu_2(\sigma).$$
		       If $\nu_2(\bar\sigma)=b_{I-1}$, then $\bar\sigma$ is one of the
differentials we
		       need to satisfy condition (1) in the statement of the Theorem and we
will
		       put $\tau_{I-1}=\bar\sigma$.  If not, then $\nu_2(\bar\sigma)=s_{k'}$
		       for some $k'$ with $k<k'\le \delta_2-1$.  Then, by adding a  suitable
		       multiple of  $\tau_{I+\delta_1+k'}$, we obtain a differential with a
greater
		       order on the second branch (while leaving the order on the first
branch
		       unchanged).  In this way, we obtain a differential, call it
$\tau_{I-1}$, such
		       that $\nu_1(\tau_{I-1})=a_{I-1}$ and $\nu_2(\tau_{I-1})=b_{I-1}$.

		       We continue by induction, assuming that we have found the
differentials
		       $\tau_{I-t+1},\dots,$ $\tau_{I-1}$ corresponding to the maximal points
		       $(a_{I-t+1},b_{I-t+1}),\dots,(a_{I-1},b_{I-1})$.  Consider the maximal
point
		       $(a_{I-t},b_{I-t})$.  One of the differentials we found above using
chain $(*)$,
		       call it $\rho$, satisfies $\nu_1(\rho)=a_{I-t}$.  If $\nu_2(\rho)\ne
b_{I-t}$, then
		       we add to $\rho$ a suitable multiple of either $\tau_{I+\delta_1+k}$
if
		       $\nu_2(\rho)=s_k$ for some $k$, or $\tau_{I-u}$ if $\nu_2(\rho)=
		       b_{I-u}$ for some $u, 1\le u\le t-1$.  Continuing in this way,
		        we can increase the order of
		       the differential on the second branch, without changing the order on
the first
		       branch, until we obtain a differential $\tau_{I-t}$ such that
$\nu_1(\tau_{I-t})=
		       a_{I-t}$ and $\nu_2(\tau_{I-t})=b_{I-t}$.  By this inductive process,
we
		       obtain differentials $\tau_0,\dots,\tau_{I-1} $ satisfying condition
(1) in the
		       statement of the Theorem.

		       The differentials $\tau_0,\tau_1,\dots,\tau_{\delta-1}$ are easily
seen to be
		       linearly independent by considering their orders on the two branches
at $P$.
		       \qed

		        \medskip
		       A basis of $g$ linearly independent dualizing differentials on $X$ may
be obtained
		       by taking the union of a basis of $g-\delta$ dualizing differentials
on $Y$ and the
		       $\delta$ differentials in Theorem (4.16).  We will divide such a
		        basis into four subsets and will use the following notation.  Let
		       $$\tau_0,\tau_1,\dots,\tau_{I-1}$$
		       denote the differentials corresponding, as in (1) of Theorem (4.16),
to the
		        maximal points of $S$.   Let
		       $$\zeta_0,\zeta_1,\dots,\zeta_{\delta_1-1}$$
		       denote the differentials corresponding, as in (2) of Theorem (4.16),
to certain
		       points in $S$ with first coordinate $\xi_1-1-l$, where $l$ is a gap of
$S_1$.
		          Note that on $Y$ each of
		        the $\zeta_j$'s is
		       regular at $Q_2$  and has a pole at $Q_1$.  Let
		       $$\eta_0,\eta_1,\dots,\eta_{\delta_2-1}$$
		       denote the differentials corresponding, as in (3) of Theorem (4.16),
to certain points
		       in $S$ with second coordinate $\xi_2-1-l'$, where $l'$ is a
		       gap of $S_2$.  On $Y$, each of the $\eta_k$'s is
		       regular at $Q_1$ and has a pole at $Q_2$. Finally, let
		       $$\sigma_0,\sigma_1,\dots,\sigma_{g-\delta-1}$$
		       be a basis of the dualizing differentials on $Y$.

		       To state the main result of this section, we must also introduce two
linear systems on
		       $Y$.  Let $$V_1\subseteq H^0(Y,\omega_Y(-c_2Q_2))$$ be the linear
system
		       generated by

$$\eta_0,,\dots,\eta_{\delta_2-1},\sigma_0,\dots,\sigma_{g-\delta-1}.$$
		       Then $V_1$ has dimension $g-\delta+\delta_2$ and
		        dim$_kH^0(Y,\omega_Y(-c_2Q_2))=g-\delta+c_2-1$, assuming
			$c_2>0$.  If $c_2=0$, then
			$V_1=H^0(Y,\omega_Y)$, while if $c_2>0$, then the codimension of
$V_1$
			in $H^0(Y,\omega_Y(-c_2Q_2))$ is $c_2-1-\delta_2$.  Hence
			 $V_1=H^0(Y,\omega_Y(-c_2Q_2))$ if and only if the semigroup
$S_2=
			\{n\in \N : n=0 \hbox{ or } n\ge c_2\}$ (e.g., if $P$ is  a simple
			cusp on the second branch).

		       Let $$V_2\subseteq H^0(Y,\omega_Y(-c_1Q_1))$$ be the linear
system
		       generated by

$$\zeta_0,\dots,\zeta_{\delta_1-1},\sigma_0,\dots,\sigma_{g-\delta-1}.$$
		       Similar remarks to those made just above also hold concerning
		        $V_2$ and $H^0(Y,\omega_Y(-c_1Q_1))$.

		       \proclaim{(4.17) Theorem}.   Suppose $X$ is a Gorenstein curve of
arithmetic
		       genus $g$.  Suppose $P$ is a singularity with precisely two branches.
		       Let $Q_1$ and $Q_2$ be the two points on the partial normalization $Y$
			of $X$ at $P$ that correspond to the branches at $P$.  Let $V_1$ and
			$V_2$ denote the linear systems on $Y$ defined above.
		       Then we have
		       $$W_X(P)=\delta(g-1)(g+1) -I(g-1)-wt(S_1)-wt(S_2)+W_{V_1}(Q_1)+
		       W_{V_2}(Q_2).$$

		       \no{\bf Proof}.  Locally at $P$, write
		       $$\eqalign{\tau_i=F_i\tau,&\qquad i=0,1,\dots,I-1\cr
		       \zeta_j=G_j\tau,&\qquad j=0,1,\dots,\delta_1-1\cr
		       \eta_k=H_k\tau,&\qquad k=0,1,\dots,\delta_2-1\cr
		       \sigma_l=M_l\tau,&\qquad l=0,1,\dots,g-\delta-1.\cr}$$
		       Put $$(\hat F,\hat G,\hat H,\hat M)=(F_0,\dots,F_{I-1},G_0,\dots,
		       G_{\delta_1-1},H_0,\dots,H_{\delta_2-1},M_0,\dots,M_{g-\delta-1}).$$
		       Then, as follows from Proposition (1.6), we have
$$W_X(P)=\delta(g-1)g
+
		        \hbox{ord}_{Q_1}W_t(\hat F,\hat G,\hat H,\hat M) +
		       \hbox{ord}_{Q_2}W_t(\hat F,\hat G,\hat H,\hat M),
		       $$
		       where $t$ is a local coordinate at $Q_1$ and $Q_2$ and $W_t$ denotes
		       the ordinary Wronskian (obtained by differentiating with respect to
$t$).
		       Notice that each of the functions $F_0,\dots,$ $F_{I-1},$
		        $G_0,\dots,G_{\delta_1-1}$
		       has a different order at $Q_1$.  By forming linear combinations of the
$H_k$'s and
		       $M_l$'s, if necessary, we may assume that each of the functions
$H_0,\dots,$
		        $H_{\delta_2-1},$ $M_0,\dots,$ $M_{g-\delta-1}$ also
			 has a different order at $Q_1$ and that,
		       of these functions, $H_0$ has the lowest order at $Q_1$, with that
order being
		       $I+2\delta_1$.
		       Then we have
		       $$\eqalign{\hbox{ord}_{Q_1}W_t(\hat F,\hat G,\hat H,\hat M)=&
		       \sum_{i=0}^{I-1}\hbox{ord}_{Q_1}F_i +
\sum_{j=0}^{\delta_1-1}\hbox{ord}_{Q_1}
		       G_j\cr
		       & + \sum_{k=0}^{\delta_2-1}\hbox{ord}_{Q_1}H_k
+\sum_{l=0}^{g-\delta-1}
		       \hbox{ord}_{Q_1}M_l -\sum_{n=0}^{g-1}n\cr
		       =&\sum_{i=0}^{I-1}(\hbox{ord}_{Q_1}F_i -i)+
		        \sum_{j=0}^{\delta_1-1}(\hbox{ord}_{Q_1}
		       G_j-(I+j))\cr
		       & + \sum_{k=0}^{\delta_2-1}(\hbox{ord}_{Q_1}H_k-(I+\delta_1+k))
		        +\sum_{l=0}^{g-\delta-1}
		       (\hbox{ord}_{Q_1}M_l -(\delta+l))\cr}$$
		       \goodbreak
		       $$\eqalign{\phantom{\hbox{ord}_{Q_1}W_t(\hat F,\hat G,\hat H,\hat
M)}
		       =&\sum_{i=0}^{I-1}(a_i-i) +
		        \sum_{j=0}^{\delta_1-1}(I+2\delta_1-1-l_{j+1}-(I+j))\cr
		       &+ \sum_{k=0}^{\delta_2-1}(\hbox{ord}_{Q_1}H_k-(I+\delta_1+k))
		        +\sum_{l=0}^{g-\delta-1}
		       (\hbox{ord}_{Q_1}M_l -(\delta+l))\cr
		       =&\> \delta_1I +(\delta_1-1)\delta_1
		        -wt(S_1)\cr
			&+\sum_{k=0}^{\delta_2-1}(\hbox{ord}_{Q_1}H_k-(I+\delta_1+k))
		        +
			\sum_{l=0}^{g-\delta-1}
		       (\hbox{ord}_{Q_1}M_l -(\delta+l)),\cr}$$
		       where $l_1,l_2,\dots,l_{\delta_1}$ are the gaps of $S_1$ and where,
		        in the last equality, we have used  Corollary (4.8) and the fact that
		       $$\sum_{j=0}^{\delta_1-1}(I+2\delta_1-1-l_{j+1}-(I+j))=
		       (\delta_1-1)\delta_1-
		       \sum_{j=1}^{\delta_1} (l_j-j).$$

		       Now, to compute $W_{V_1}(Q_1)$, we must express
$\eta_0,\dots,\eta_{\delta_2-1},
		       \sigma_0,\dots,\sigma_{g-\delta-1}$ in terms of a generator of
$\omega_Y(-c_2Q_2)$
		       at $Q_1$.
		      Let $\eta$ be a generator of $\omega_Y(-c_2Q_2)$
		      at $Q_1$.  Then $\eta$ has order 0 at $Q_1$.  Note that the rational
function
		       $H=\eta/\tau$ has
		      a zero of order $I+2\delta_1$ at $Q_1$.  We then have

$$\eqalign{W_{V_1}(Q_1)=&\hbox{ord}_{Q_1}W_t(H_0/H,\dots,H_{\delta_2-1}/H,
		      M_0/H,\dots,M_{g-\delta-1}/H)\cr
=&\sum_{k=0}^{\delta_2-1}((\hbox{ord}_{Q_1}H_k-
(I+2\delta_1))-k)+\sum_{l=0}^
		      {g-\delta-1}((\hbox{ord}_{Q_1}M_l-(I+2\delta_1))-(\delta_2+l))\cr
=&\sum_{k=0}^{\delta_2-1}(\hbox{ord}_{Q_1}H_k-
(I+\delta_1+k))-\delta_2\delta_1\cr
			&+\sum_{l=0}^{g-\delta-1}(\hbox{ord}_{Q_1}M_l-(\delta+l)) -
(g-\delta)\delta_1.\cr}$$
			Thus,
			$$\eqalign{\hbox{ord}_{Q_1}&W_t(\hat F,\hat G,\hat H,\hat M)\cr
			&=\delta_1I+
(\delta_1-1)\delta_1-wt(S_1)+\delta_2\delta_1+
(g-\delta)\delta_1+W_{V_1}(Q_1)\cr
			&=\delta_1(I+\delta_1-1+\delta_2+g-\delta)-wt(S_1)+W_{V_1}(Q_1)
\cr
			&=\delta_1(g-1)-wt(S_1)+W_{V_1}(Q_1),\cr}$$
			since $\delta=I+\delta_1+\delta_2$.

			Similarly, we have
			$$\hbox{ord}_{Q_2}W_t(\hat F,\hat G,\hat H,\hat
M)=\delta_2(g-1)-wt(S_2)+
			W_{V_2}(Q_2).$$
			The Theorem now follows by adding these two orders and using the
fact that
			$\delta_1+\delta_2=\delta-I$.

			\line{\hfil\qed}
			\medskip
			In the case of an ordinary node, we have $I=1,\delta_1=\delta_2=0$
and
Theorem
			(4.17) reduces to the following result of Widland [21].

			\proclaim {(4.18) Corollary}. If $P$ is an ordinary node, then
			$$W_X(P)=(g-1)g+W_Y(Q_1)+W_Y(Q_2).$$

			\medskip
			We will call a singularity $P$ {\it overweight\/} if its Weierstrass
weight
			 is greater than
			the ``expected'' number.  A unibranch singularity  $P$ is not
overweight if
			the point lying over $P$ on the partial normalization at $P$ is not a
Weierstrass
			point (see Theorem (3.4)).
			 A singularity $P$ with two branches is not overweight if, with the
notation
			of Theorem (4.17), $W_{V_1}(Q_1) =W_{V_2}(Q_2)=0.$  We can
now state a
			result in the case of a two-branch singularity that is analogous to (3) of
			Theorem (3.5).

			\proclaim{(4.19) Proposition}. Suppose  $X$ is a rational Gorenstein
curve
			of arithmetic genus $g$ with a
			single singularity $P$.  Suppose that $P$ has precisely two
			 branches and let $I, S_1$, and
			$S_2$ be as defined above.   If $P$ is not overweight,  then the number
of
smooth
			Weierstrass points on $X$, counting multiplicities, is
			$$I(g-1)+wt(S_1) + wt(S_2).$$

			\no{\bf Proof.}   Since $X$ is rational, we have $g=\delta$.  The weight
of
$P$ is
			then $g^3-g-I(g-1)-wt(S_1)-wt(S_2)$, as follows from Theorem (4.17)
			since $P$ is not overweight.\qed
			\medskip
			One can also prove a result similar to Theorem (3.7) in the case of
			a rational curve with unibranch and  two-branch singularities that are
			 not overweight.  The result is that the number of smooth Weierstrass
			 points, counting multiplicities, is given as the sum of local
contributions
			 from the singularities, with each unibranch singularity contributing the
			 weight of its semigroup and each two-branch singularity contributing
			 $I(g-1)+wt(S_1)+wt(S_2)$.

			\bigskip
			\goodbreak
			\centerline{\bf References}
			\bigskip
			\frenchspacing
			\item{1.} Delgado, F.: The semigroup of values of a curve singularity
			with several branches, Manuscripta Math. {\bf 59}, 347--374 (1987)
			\item{2.} Delgado,F.: Gorenstein  curves and symmetry
			of the semigroup of values, Manuscripta Math. {\bf 61}, 285--296
(1988)
			\item{3.}  Farkas, H., Kra, I.: Riemann surfaces,  New York: Springer
			1980
			\item{4.} Garcia, A.: Semigroups associated to singular points of plane
curves,
			J. Reine Angew. Math. {\bf 336}, 165--184 (1982)
			\item{5.} Hartshorne, R.: Generalized divisors on Gorenstein
			curves and a theorem of Noether, J. Math. Kyoto Univ. {\bf 26},
375--386
			(1986)
			\item{6.}  Henn, H.-W.: Minimalanzahl von Weierstrasspunkten,
Archiv
			der Math. {\bf 31}, 38--43 (1978)
			\item{7.} Hironaka, H.: On the arithmetic genera and the effective
genera
			of algebraic curves, Memoirs of the College of Science, Univ. Kyoto,
			Series A, {\bf 30}, 178--195 (1957)
			\item{8.}  Homma, M. Duality of space curves and their tangent
surfaces
			in characteristic $p>0$, Ark. Mat. {\bf 29}, 221--235 (1991)
			\item{9.} Kunz, E: The value-semigroup of a one-dimensional
Gorenstein
			ring, Proc. Amer. Math. Soc. {\bf 25}, 748--751 (1970)
			\item{10.} Laksov, D.: Weierstrass points on curves, Ast\'erisque
			 {\bf 87--88}, 221--247 (1981)
			 \item{11.} Laksov, D., Thorup, A.: The Brill-Segre formula for
families of
curves.
			 In: Enumerative algebraic geometry (Copenhagen, 1989), 131--148,
Contemp.
			 Math. {\bf 123}. Providence, RI: Amer. Math. Soc. 1991
			\item{12.}  Laksov, D. , Thorup, A.:  On gap sequences for families of
curves,
			preprint
			\item{13.} Lax, R.F.: On the distribution of Weierstrass points on
singular
			curves,  Israel J. Math. {\bf 57}, 107--115 (1987)
			\item{14.} Lax, R.F., Widland, C.: Gap sequences at a singularity, Pac.
J.
			Math. {\bf 150}, 111--122 (1991)
			\item{15.} Oliveira, G.: Weierstrass semigroups and the canonical
ideal
			of non-trigonal curves, Manuscripta Math. {\bf 71}, 431--450 (1991)
			\item{16.} Schmidt, F.K.:  Die Wronskiche Determinante in beliebigen
			differenzierbaren Funktionenk\"orpern, Math. Z. {\bf 45}, 62--74 (1939)
			\item{17.} Schmidt, F.K.: Zur arithmetischen Theorie der
			algebraischen Funktionen II.  Allgemeine Theorie der
Weierstrasspunkte,
			Math. Z. {\bf 45}, 75--96 (1939)
			\item{18.} Serre, J.-P.: Groupes Alg\'ebriques et Corps des Classes.
Paris:
			 Hermann 1959
			\item{19.} St\"ohr, K.-O., Voloch, J.F.: Weierstrass points and curves
over
			finite fields, Proc. London Math. Soc. (3) {\bf 52}, 1--19 (1986)
			\item{20.} St\"ohr, K.-O.: On the moduli spaces of Gorenstein curves
with
			symmetric semigroups, preprint
			\item{21.} Widland, C.: On Weierstrass points of Gorenstein curves,
Ph.D.
			dissertation, Louisiana State University 1984
			\item{22.} Widland, C.,  Lax, R.F.: Weierstrass points on Gorenstein
curves,
			Pac. J. Math. {\bf 142}, 197--208 (1990)

			\bye